\documentclass{article}

\usepackage[a4paper]{geometry}
\usepackage[indent]{parskip}
\usepackage{etoolbox}
\usepackage{ifdraft}

\usepackage[dvipsnames]{xcolor}
\definecolor{Basil}{HTML}{32612D}
\definecolor{Berry}{HTML}{7A1712}

\usepackage[author={}]{fixme}
\fxsetup{marginface=\scriptsize\sffamily\color{gray}}

\usepackage[mathlines]{lineno}
\ifdraft{\linenumbers}{}

\usepackage[british]{babel}
\usepackage[style=american]{csquotes}

\usepackage{lmodern}
\usepackage[T1]{fontenc}
\usepackage[sf,bf,nobottomtitles*,]{titlesec}
\usepackage{abstract}
\titleformat{\section}{\normalfont\large\bfseries\sffamily}{\thesection}{1em}{}
\titleformat{\subsection}{\normalfont\normalsize\bfseries\sffamily}{\thesubsection}{1em}{}

\usepackage[nodisplayskipstretch]{setspace}
\usepackage{microtype}
\usepackage[prevent-all]{widows-and-orphans}
\usepackage{titling}
\pretitle{\begin{center}\begin{minipage}{0.8\textwidth}\Large\centering}
\posttitle{\end{minipage}\end{center}\vskip 1em}
\preauthor{\begin{center}\large \lineskip 0.5em\begin{tabular}[t]{c}}
\postauthor{\end{tabular}\par\end{center}\vskip 1em}
\usepackage{cite}  
\newcounter{cite}
\preto{\cite}{\stepcounter{cite}}{}{}
\NewCommandCopy{\oldbibliography}{\bibliography}
\renewcommand{\bibliography}[1]{\ifnum\value{cite}>0\oldbibliography{#1}\fi}

\newcommand{\keywords}[1]{
    \begingroup
    \def\and{ ; }
    \hypersetup{pdfkeywords={#1}}
    \def\and{\ifhmode\unskip\nobreak\fi\ $\cdot$ }
    \paragraph*{Keywords} #1
    \endgroup
}
\newcommand{\msc}[1]{
    \begingroup
    \def\and{ ; }
    \hypersetup{pdfsubject={#1}}
    \def\and{\ifhmode\unskip\nobreak\fi\ $\cdot$ }
    \paragraph*{Mathematics Subject Classification (2010)} #1
    \endgroup
}

\usepackage{mathtools}
\usepackage{amssymb}
\usepackage{cases}
\usepackage[bb=dsfontserif]{mathalpha}

\usepackage{amsthm}
\theoremstyle{definition}

\newtheorem{definition}{Definition}[section]

\theoremstyle{plain}

\newtheorem{proposition}{Proposition}[section]
\newtheorem{lemma}{Lemma}[section]
\newtheorem{theorem}{Theorem}[section]

\theoremstyle{remark}

\mathcode`@="8000{\catcode`\@=\active\gdef@{\mkern1mu}}

\usepackage[pdfusetitle,hyperfootnotes=false]{hyperref}
\usepackage{url}
\hypersetup{
    final,
    colorlinks=true,
    linkcolor=Berry,
    citecolor=Basil,
    urlcolor=RedViolet,
}
\newcommand*{\email}[2][inbox]{Email: \href{mailto:#2}{\texttt{#2}}}
\newcommand*{\orcid}[1]{ORCID: \href{https://orcid.org/#1}{\texttt{#1}}}

\usepackage[acronym]{glossaries}
\glsdisablehyper
\newacronym{dfo}{DFO}{derivative-free optimization}

\DeclareMathOperator{\bigo}{\mathcal{O}}
\DeclareMathOperator{\card}{card}
\newcommand*{\abs}[2][]{#1\lvert#2#1\rvert}
\newcommand*{\norm}[2][]{#1\lVert#2#1\rVert}
\newcommand*{\set}[2][]{#1\{#2#1\}}
\newcommand*{\obj}{f}
\newcommand*{\objm}[1][]{\hat{\obj}\ifblank{#1}{}{^{#1}}}
\newcommand*{\fset}{\Omega}
\newcommand*{\R}{\mathbb{R}}
\newcommand*{\T}{\mathsf{T}}
\newcommand*{\lagp}[1][]{L\ifblank{#1}{}{_{#1}}}
\newcommand*{\xpt}[1][]{\mathcal{Y}\ifblank{#1}{}{^{@#1}}}
\newcommand*{\ones}{\mathbb{1}}
\newcommand*{\qpoly}{\mathcal{Q}_n}

\title{An Optimal Interpolation Set for Model-Based Derivative-Free Optimization Methods}
\author{
    Tom M. Ragonneau%
    \thanks{Department of Applied Mathematics, The Hong Kong Polytechnic University, Hong Kong, China.}%
    \thanksgap{1ex}%
    \thanks{\email{tom.ragonneau@gmail.com}; \orcid{0000-0003-2717-2876}.}
    \and
    Zaikun Zhang%
    \texorpdfstring{\thanksmark{1}}{}%
    \thanksgap{1ex}%
    \thanks{\email{zaikunzhang@gmail.com}; \orcid{0000-0001-8934-8190}.}
}
\date{\today}

\begin{document}

\maketitle

\begin{abstract}
    This paper demonstrates the optimality of an interpolation set employed in derivative-free trust-region methods.
    This set is optimal in the sense that it minimizes the constant of well-poisedness in a ball centred at the starting point.
    It is chosen as the default initial interpolation set by many derivative-free trust-region methods based on underdetermined quadratic interpolation, including NEWUOA, BOBYQA, LINCOA, and COBYQA.
    Our analysis provides a theoretical justification for this choice.
\end{abstract}


\msc{90C56 \and 41A10 \and 65K05 \and 90C30}

\section{Introduction}

\Gls{dfo} methods solve the optimization problem
\begin{equation}
    \label{eq:nlp}
    \min_{x \in \fset \subseteq \R^n} @@ \obj(x)\\
\end{equation}
based on function evaluations without relying on derivative information.
They are needed when classical or generalized derivatives of the objective function~$\obj \colon \R^n \to \R$ are unavailable or too expensive to evaluate, which may also be the case for the constraint functions underlying the feasible set~$\fset \subseteq \R^n$.
Examples of such problems arise from many areas, such as reinforcement learning~\cite{Qian_Yu_2021}, hyperparameter tuning~\cite{Ghanbari_Scheinberg_2017}, particle physics~\cite{Eldred_Etal_2023}, and aircraft engineering~\cite{Gazaix_Etal_2019}.
The \gls{dfo} literature is vast.
For comprehensive overviews, we refer to the monographs~\cite{Conn_Scheinberg_Vicente_2009,Audet_Hare_2017}, the survey papers~\cite{Powell_1975,Powell_1998,Rios_Sahinidis_2013,Custodio_Scheinberg_Vicente_2017,Larson_Menickelly_Wild_2019}, the recent thesis~\cite{Ragonneau_2022}, and the references therein.

This paper studies a set employed as the initial interpolation set by many derivative-free trust-region methods~\cite{Conn_Gould_Toint_2000,Yuan_2015} based on underdetermined quadratic interpolation, such as NEWUOA~\cite{Powell_2006}, and, when no constraints are present, BOBYQA~\cite{Powell_2009}, LINCOA, and COBYQA~\cite{Ragonneau_2022,Ragonneau_Zhang_2023}.
It is also employed by~\cite{Conejo_Karas_Pedroso_2015} in a method for constrained problems.
We show that this interpolation set is optimal in the sense that it minimizes the constant of well-poisedness in a ball centred at the starting point, which is detailed in Theorem~\ref{thm:optimset}.

In what follows, we denote by~$\qpoly$ the space of polynomials on~$\R^n$ of degree at most two, which will be referred to as quadratic polynomials for simplicity.
The \gls{dfo} methods that we consider maintain an interpolation set~$\xpt[k] \subseteq \R^n$ at iteration~$k$, and build a model~$\objm[k] \in \qpoly$ for~$\obj$ according to the interpolation conditions
\begin{equation}
    \label{eq:intp-cond}
    \objm[k](y) = \obj(y), \quad y \in \xpt[k].
\end{equation}
Clearly, if the interpolation conditions~\eqref{eq:intp-cond} do not contradict each other, the model~$\objm[k]$ is uniquely defined by~\eqref{eq:intp-cond} only if
\begin{equation*}
    \card(\xpt[k]) = \dim(\qpoly) = \frac{1}{2} (n + 1) (n + 2).
\end{equation*}
Such a scheme, referred to as the fully-determined quadratic interpolation, is used, for example, by UOBYQA~\cite{Powell_2002}.
However, the first iteration of this method requires~$\bigo(n^2)$ function evaluations to construct the initial model~$\objm[@0]$.
This is impracticable unless~$n$ is small.
Therefore, modern interpolation-based \gls{dfo} methods work with fewer interpolation points, giving rise to the underdetermined quadratic interpolation.
To uniquely define~$\objm[k]$ in this case, it is common to solve an interpolation problem of the form
\begin{equation}
    \label{eq:intp-var}
    \begin{aligned}
        \min_{Q \in \qpoly} & \quad \mathcal{F}^{@k}(Q)\\
        \text{s.t.}         & \quad Q(y) = \obj(y), \quad y \in \xpt[k],
    \end{aligned}
\end{equation}
where~$\mathcal{F}^{@k}$ is a functional that promotes certain desired regularity of~$\objm[k]$.
For example, DFO~\cite{Conn_Scheinberg_Toint_1998} and MNH~\cite{Wild_2008} take
\begin{equation}
    \label{eq:mnh}
    \mathcal{F}^{@k}(Q) = \norm[\big]{\nabla^2 Q}_{\mathsf{F}}^2,
\end{equation}
where~$\norm{\cdot}_{\mathsf{F}}$ denotes the Frobenius norm.
In contrast, NEWUOA, BOBYQA, LINCOA, and COBYQA employ
\begin{equation}
    \label{eq:df-psb}
    \mathcal{F}^{@k}(Q) = \norm[\big]{\nabla^2 Q - \nabla^2 \objm[k - 1]}_{\mathsf{F}}^2
\end{equation}
with~$\objm[-1] = 0$.
There are other functionals based on the~$\ell_1$-norm~\cite{Bandeira_Scheinberg_Vicente_2012} or Sobolev seminorms of quadratic polynomials~\cite{Zhang_2014,Xie_Yuan_2022}.
The variational problem described by~\eqref{eq:intp-var} and~\eqref{eq:df-psb} defines the derivative-free symmetric
Broyden update proposed by Powell (see~\cite{Powell_2004a,Powell_2013},~\cite[\S~3.6]{Zhang_2012}, and~\cite[\S~2.4.2]{Ragonneau_2022}), with the functional~\eqref{eq:df-psb} being inspired by the least-change property of quasi-Newton updates~\cite{Dennis_Schnabel_1979}.
To ensure that~\eqref{eq:intp-var} admits a solution that includes second-order information, one normally requires that
\begin{equation*}
    n + 2 \le \card(\xpt[k]) \le \frac{1}{2} (n + 1) (n + 2).
\end{equation*}
Finally, if the constraint functions underlying the feasible set~$\fset$ in~\eqref{eq:nlp} are nonlinear, we can approximate them using similar techniques, which is the case for COBYQA.

Most interpolation-based \gls{dfo} methods do not build the interpolation point~$\xpt[k]$ \emph{ab initio} at each iteration.
Instead, they choose an initial interpolation set~$\xpt[0]$ and then  update~$\xpt[k]$ to obtain~$\xpt[k + 1]$.
For example, only one point differs between~$\xpt[k]$ and~$\xpt[k + 1]$ in Powell's trust-region \gls{dfo} methods and COBYQA during usual iterations.
This mechanism tries to make the best of the function evaluations, which are normally considered the principal cost of \gls{dfo} methods.
In such a case, it is clear that the initial interpolation set~$\xpt[0]$ must be carefully chosen, as it will impact the optimization method for subsequent iterations.
Powell devised an initial interpolation set for his methods, described in Subsection~\ref{subsec:powell-set}.
Although the design of this set is natural, no theoretical analysis has been provided to justify its configuration.
This paper contributes such an analysis, showing that no other interpolation set is better in a well-poisedness sense~\cite[\S~3.3]{Conn_Scheinberg_Vicente_2009}.

This paper is organized as follows.
Section~\ref{sec:well-poisedness} introduces the $\Lambda$-poisedness of interpolation sets in the minimum Frobenius norm sense.
Section~\ref{sec:main-result} then employs this notion to study the initial interpolation set used by Powell.
We provide lower and upper bounds for the constant of well-poisedness of the set, evaluate this constant in a few special cases, and then demonstrate the optimality of this set under the default setting recommended by Powell.
In addition, we raise some open questions motivated by our results.
Section~\ref{sec:conclusion} concludes this paper.

\section{Well-poisedness of interpolation sets}
\label{sec:well-poisedness}

In this section, we consider an interpolation set
\begin{equation*}
    \xpt = \set{y^1, @@ y^2, @@ \dots, @@ y^m} \subseteq \R^n
\end{equation*}
and the interpolation problem
\begin{equation}
    \label{eq:intp-cond-gen}
    \begin{aligned}
        \min_{Q \in \qpoly} & \quad \norm{\nabla^2 Q}_\textsf{F}\\
        \text{s.t.}         & \quad Q(y) = \obj(y), \quad y \in \xpt.
    \end{aligned}
\end{equation}
This problem covers~\eqref{eq:intp-var}--\eqref{eq:mnh}, and also covers~\eqref{eq:intp-var} with~\eqref{eq:df-psb} by a simple change of variable.

\begin{definition}[Poisedness]
    The set~$\xpt$ is \emph{poised} in the minimum Frobenius norm sense if problem~\eqref{eq:intp-cond-gen} has a unique solution for any real-valued function~$\obj$.
\end{definition}

Observe that the variational problem~\eqref{eq:intp-cond-gen} is an equality-constrained quadratic programming problem with respect to the coefficients of~$Q$.
Thus, its KKT system is linear.
The exact formulation of this system can be found in~\cite{Powell_2004a,Powell_2004b}.
The interpolation set~$\xpt$ is poised if this KKT system is uniquely solvable (see~\cite[\S~2]{Powell_2004a} and~\cite[\S~5.3]{Conn_Scheinberg_Vicente_2009}).
Intuitively,~$\xpt$ can be said well-poised if this system is well-conditioned.
This intuition will be formalized in the sequel.

\subsection{Minimum Frobenius norm Lagrange polynomials}

To formally define a measure of well-poisedness of~$\xpt$, we first need to extend the classical definition of the Lagrange polynomials as follows.

\begin{definition}[Minimum Frobenius norm Lagrange polynomials~{\cite[Definition~5.1]{Conn_Scheinberg_Vicente_2009}}]
    \label{def:min-norm-lagp}
    Assume that the interpolation set~$\xpt$ is poised.
    The~$i$th \emph{minimum Frobenius norm Lagrange polynomial}~$\lagp[i]$ for the interpolation problem~\eqref{eq:intp-cond-gen}, with~$i \in \set{1, @@ 2, @@ \dots, @@ m}$, is the unique quadratic polynomial that solves
    \begin{equation*}
        \begin{aligned}
            \min_{Q \in \qpoly} & \quad \norm[\big]{\nabla^2 Q}_{\mathsf{F}}\\
            \text{s.t.}         & \quad Q(y^j) = \delta_{i, j}, \quad j \in \set{1, @@ 2, @@ \dots, @@ m},
        \end{aligned}
    \end{equation*}
    where~$\delta_{i, j}$ denotes the Kronecker delta.
\end{definition}

Remark that if~$m = (n + 1) (n + 2) / 2$, then Definition~\ref{def:min-norm-lagp} reduces to the classical definition of the Lagrange polynomials.
Moreover, it can be shown that the solution~$\objm$ to problem~\eqref{eq:intp-cond-gen} is given by
\begin{equation*}
    \objm(x) = \sum_{i = 1}^m \obj(y^i) \lagp[i](x)
\end{equation*}
for any~$x \in \R^n$~\cite[Lemma~5.2]{Conn_Scheinberg_Vicente_2009}.

\subsection{Well-poisedness in the minimum Frobenius norm sense}

We are now equipped to define the notion of~$\Lambda$-poisedness in the minimum Frobenius norm sense.

\begin{definition}[$\Lambda$-poisedness in the minimum Frobenius norm sense~{\cite[Definition~5.6]{Conn_Scheinberg_Vicente_2009}}]
    \label{def:lambda-p}
    Let~$\Lambda$ be a positive constant and~$\mathcal{C} \subseteq \R^n$ be a compact set.
    The set~$\xpt$ is said to be~\emph{$\Lambda$-poised in the minimum Frobenius norm sense} in~$\mathcal{C}$ if it is poised and
    \begin{equation*}
        \max_{1 \le i \le m} \max_{x \in \mathcal{C}} @@ \abs[\big]{\lagp[i](x)} \le \Lambda.
    \end{equation*}
\end{definition}

One can show that~$\xpt$ is~$\Lambda$-poised in the unit ball if and only if the condition number of the KKT system of~\eqref{eq:intp-var} is bounded by some polynomials of~$\Lambda$~\cite[Theorem~5.8]{Conn_Scheinberg_Vicente_2009}.
Therefore, the notion of~$\Lambda$-poisedness formalizes the intuition about well-poisedness mentioned above.

If~$\xpt$ is a poised interpolation set and~$\mathcal{C} \subseteq \R^n$ is compact, then we will refer to
\begin{equation*}
    \Lambda_{\mathcal{C}}(\xpt) = \max_{1 \le i \le m} \max_{x \in \mathcal{C}} @@ \abs[\big]{\lagp[i](x)}
\end{equation*}
as the \emph{constant of well-poisedness} of~$\xpt$ in~$\mathcal{C}$.

\section{Optimality of Powell's initial interpolation set}
\label{sec:main-result}

This section first presents the initial interpolation set devised by Powell for his trust-region \gls{dfo} methods, and then analyzes its well-poisedness.
We establish bounds on its constant of well-poisedness and evaluate this constant in some special cases.
Finally, we point out that the default setting in Powell's methods renders an optimal interpolation set in terms of the well-poisedness in a ball centred at the starting point.

\subsection{Description of the initial interpolation set}
\label{subsec:powell-set}

The initial interpolation that Powell designed for his methods is as follows.
Without loss of generality, we assume that the starting point is at the origin.
Suppose that~$\Delta > 0$ is the initial trust-region radius.
For~$i \in \set{1, @@ 2, @@ \dots, @@ 2n + 1}$, define
\begin{equation}
    \label{eq:set-def}
    y^i =
    \begin{cases}
        0,                      & \text{if~$i = 1$,}\\
        \Delta e_{i - 1},       & \text{if~$2 \le i \le n + 1$,}\\
        -\Delta e_{i - n - 1},  & \text{otherwise},
    \end{cases}
\end{equation}
where~$e_i \in \R^n$ denotes the $i$th canonical coordinate vector in~$\R^n$.
Let~$m$ be the number of interpolation points.
We focus on the case with~$n + 2 \le m \le 2n + 1$ following Powell's suggestion.\footnote{Powell's BOBYQA code contains a comment that \textquote[][.]{choices that exceed $2n+1$ are not recommended}}
The initial interpolation set is then chosen to be
\begin{equation}
    \label{eq:set-def-m}
    \xpt[0]_m = \set{y^1, @@ y^2, @@ \dots, @@ y^m},
\end{equation}
which can be found in~\cite[Equation~(3.2)]{Powell_2006}.
The default value for~$m$ proposed by Powell is~$2n + 1$.
In other words, the default initial interpolation set is~$\xpt[0] = \xpt[0]_{2n + 1}$.
This is a natural choice, as the corresponding interpolation set is geometrically appealing.
It consists of the origin and equidistant points to the origin in the positive and negative coordinate directions.
In what follows, we provide a theory that justifies this natural choice, showing that~$\xpt[0]_{2n + 1}$ is indeed optimal in terms of well-poisedness.

\subsection{Well-poisedness of the interpolation set}

We now investigate the $\Lambda$-poisedness of~$\xpt[0]_m$.
We will do this in the closed~$\ell_p$-norm ball of radius~$\Delta$ centred at the origin, namely
\begin{equation*}
    \mathcal{B}_p(\Delta) = \set{x \in \R^n : \norm{x}_p \le \Delta},
\end{equation*}
with~$p \in [1, \infty]$.
Note that we allow~$p = \infty$, and that~$\mathcal{B}_p(\Delta)$ is also the smallest~$\ell_p$-norm ball enclosing~$\xpt[0]_m$.
Trust-region \gls{dfo} methods usually define the trust region by the Euclidean norm, so the case with~$p=2$ is the most
interesting.
However, it can be beneficial to define the trust region by polyhedral norms when bound or linear constraints are present, and $p = 1$ or~$p = \infty$ will become more relevant.

According to Definition~\ref{def:lambda-p}, the set~$\xpt[0]_m$ is~$\Lambda_p$-poised in~$\mathcal{B}_p(\Delta)$ in the minimum Frobenius norm sense with
\begin{equation}
    \label{eq:Lambdap}
    \Lambda_p = \max_{1 \le i \le m} \max_{\norm{x}_p \le \Delta} @@ \abs[\big]{\lagp[i](x)},
\end{equation}
where~$\lagp[i]$, for~$i \in \set{1, @@ 2, @@ \dots, @@ m}$, is the~$i$th minimum Frobenius norm Lagrange polynomial associated with~$\xpt[0]_m$.
Indeed, $\Lambda_p$ is the {constant of well-poisedness} of~$\xpt[0]_m$ in~$\mathcal{B}_p(\Delta)$.
We will focus on this constant in what follows.

\subsubsection{Formulation of the Lagrange polynomials}

To study~$\Lambda_p$, we first present explicit formulae for~$\lagp[i]$ for all~$i \in \set{1, @@ 2, @@ \dots, @@ m}$.
These formulae are given in~\cite[\S~3]{Powell_2006} without proof.

\begin{lemma}
    \label{lem:lagp}
    For all~$x \in \R^n$ and~$m \in \set{n + 2, @@ n + 3, @@ \dots,  @@ 2n + 1}$, the expression of~$\lagp[i]$, for~$i \in \set{1, @@ 2, @@ \dots, @@ m}$, is given by
    \begin{equation}
        \label{eq:lagp}
        \lagp[i](x) =
        \begin{cases}
            \displaystyle 1 - \frac{1}{\Delta^2} \sum_{j = 1}^{m - n - 1} \!\! x_j^2 - \frac{1}{\Delta} \sum_{j = m - n}^n \!\! x_j,    & \text{if~$i = 1$,}\\[2ex]
            \displaystyle \frac{x_{i - 1}^2}{2\Delta^2} + \frac{x_{i - 1}}{2\Delta},                                                    & \text{if~$2 \le i \le m - n$,}\\[2ex]
            \displaystyle \frac{x_{i - 1}}{\Delta},                                                                          & \text{if~$m - n + 1 \le i \le n + 1$,}\\[2ex]
            \displaystyle \frac{x_{i - n - 1}^2}{2\Delta^2} - \frac{x_{i -n- 1}}{2\Delta},                                    & \text{if~$n+2\le i \le m$.}
        \end{cases}
    \end{equation}
    Here, $x_j$ denotes the~$j$th entry of~$x$ for each~$j \in \set{1, @@ 2, @@ \dots, @@ n}$, and we define~$\sum_{j = m - n}^n x_j = 0$ in the formulation of~$\lagp[1]$ if~$m = 2n + 1$.
\end{lemma}

\begin{proof}
    Let~$i \in \set{1, @@ 2, @@ \dots, @@ m}$ be fixed and let~$\lagp$ be a quadratic polynomial satisfying
    \begin{equation}
        \label{eq:lagp-p}
        \lagp(y^j) = \delta_{i, j}, \quad j \in \set{1, @@ 2, @@ \dots, @@ m}.
    \end{equation}
    First, it is straightforward to verify that~$\lagp[i]$ satisfies the interpolation conditions~\eqref{eq:lagp-p}.
    Hence, it suffices to show that~$\norm{\nabla^2 \lagp[i]}_\mathsf{F} \le \norm{\nabla^2 \lagp}_\mathsf{F}$.

    Consider any~$j \in \set{1, @@ 2, @@ \dots, @@ m - n - 1}$.
    Denote the $j$th diagonal entries of~$\nabla^2 \lagp$ and~$\nabla^2 \lagp[i]$ by~$(\nabla^2 \lagp)_{j, j}$ and~$(\nabla^2 \lagp[i])_{j, j}$, respectively.
    According to equation~\eqref{eq:set-def}, we have
    \begin{equation*}
        y^1 = 0, \quad y^{j + 1} = \Delta e_j, \quad \text{and} \quad y^{n + j + 1} = -\Delta e_j.
    \end{equation*}
    Since~$\lagp$ and~$\lagp[i]$ are quadratic polynomials sharing the same values on~$\set{y^1, @@ y^{j + 1}, @@ y^{n + j + 1}}$, Taylor expansions of the quadratic polynomial~$\lagp$ around~$y^1$ yield
    \begin{equation*}
        \big(\nabla^2 \lagp \big)_{j, j} = \frac{\lagp(y^{j + 1}) + \lagp(y^{n + j + 1}) - 2 \lagp(y^1)}{\Delta^2} = \big(\nabla^2 \lagp[i] \big)_{j, j}.
    \end{equation*}
    On the other hand, it is easy to check according to~\eqref{eq:lagp} that all the entries of~$\nabla^2 \lagp[i]$ are zero except for the first~$m-n-1$ diagonal entries.
    Therefore,
    \begin{equation*}
        \norm[\big]{\nabla^2 \lagp[i]}_{\mathsf{F}}^2 \le \norm[\big]{\nabla^2 \lagp}_{\mathsf{F}}^2,
    \end{equation*}
    which completes the proof.
\end{proof}

\subsubsection{Bounds for the constant of well-poisedness}

The next lemma simplifies the expression of~$\Lambda_p$ defined in~\eqref{eq:Lambdap} for further computations.

\begin{lemma}
    \label{lem:lambda-p}
    For any~$m \in \set{n + 2, @@ n + 3, @@ \dots, @@ 2n + 1}$ and any~$p \in [1, \infty]$, we have
    \begin{equation}
        \label{eq:lambda-p}
        \Lambda_p = \max_{\norm{x}_p \le \Delta} @@ \abs[\big]{\lagp[1](x)}.
    \end{equation}
\end{lemma}

\begin{proof}
    According to Lemma~\ref{lem:lagp}, for each~$i \in \set{2, @@ 3, @@ \dots, @@ n + 1}$,~$\lagp[i](x)$ only depends on~$x_{i - 1}$ for all~$x \in \R^n$, and hence
    \begin{equation*}
        \max_{\norm{x}_p \le \Delta} @@ \abs[\big]{\lagp[i](x)} = \max_{t \in [-\Delta, \Delta]} @@ \abs[\big]{\lagp[i](t e_{i - 1})} = 1.
    \end{equation*}
    Similarly, for each~$i \in \set{n + 2, @@ n + 3, @@ \dots, @@ m}$, since~$\lagp[i](x)$ only depends on~$x_{i - n - 1}$ for all~$x \in \R^n$, we have
    \begin{equation*}
        \max_{\norm{x}_p \le \Delta} @@ \abs[\big]{\lagp[i](x)} = \max_{t \in [-\Delta, \Delta]} @@ \abs[\big]{\lagp[i](t e_{i - n - 1})} = 1.
    \end{equation*}
    Meanwhile, since~$\lagp[1](y^1) = 1$ and~$y^1 \in \mathcal{B}_p(\Delta)$, we have
    \begin{equation*}
        \max_{\norm{x}_p \le \Delta} @@ \abs[\big]{\lagp[1](x)} \ge \lagp[1](y^1) = 1.
    \end{equation*}
    Hence,~\eqref{eq:lambda-p} holds according to the definition of~$\Lambda_p$ in~\eqref{eq:Lambdap}.
\end{proof}

We are now equipped to develop lower and upper bounds for~$\Lambda_p$.
For convenience, we define henceforth
\begin{equation*}
    0^0 = 0 \quad \text{and} \quad \frac{\infty}{\infty} = 1.
\end{equation*}

\begin{theorem}
    \label{thm:bounds}
    For any~$m \in \set{n + 2, @@ n + 3, @@ \dots, @@ 2n + 1}$ and any~$p \in [1, \infty]$, we have
    \begin{equation}
        \label{eq:bounds}
        1 + (2n + 1 - m)^{\frac{p - 1}{p}} \le \Lambda_p \le n.
    \end{equation}
    In particular, we have~$\Lambda_{\infty} = \max @ \set{n - 1, @@ 2n - m + 2}$.
\end{theorem}

\begin{proof}
    We will establish the bounds in~\eqref{eq:bounds} using the formulation of~$\Lambda_p$ in Lemma~\ref{lem:lambda-p}.
    For the lower bound, by considering only the points in~$\R^n$ whose leading~$m - n - 1$ entries are zeros and whose remaining~$2n + 1 - m$ entries are equal, we have
    \begin{equation*}
        \Lambda_p = \max_{\norm{x}_p \le \Delta} @@ \abs{\lagp[1](x)} \ge \max_{t \in \R} \set[\bigg]{1 - \frac{1}{\Delta} (2n + 1 - m) t : (2n + 1 - m)^{\frac{1}{p}}\abs{t} \le \Delta} = 1 + (2n + 1 - m)^{\frac{p - 1}{p}}.
    \end{equation*}
    We now establish the upper bound.
    For any~$p \ge 1$, we have~$\mathcal{B}_p(\Delta) \subseteq \mathcal{B}_{\infty}(\Delta)$, so that~$\Lambda_p \le \Lambda_{\infty}$.
    Therefore, we only need to show that~$\Lambda_{\infty} \le \max @ \set{n - 1, @@ 2n - m + 2} \le n$.
    For any~$x \in \mathcal{B}_{\infty}(\Delta)$,
    \begin{equation}
        \label{eq:bounds-proof-1}
        -\lagp[1](x) = -1 + \frac{1}{\Delta^2} \sum_{j = 1}^{m - n - 1} \!\! x_j^2 + \frac{1}{\Delta} \sum_{j = m - n}^n \!\! x_j \le n-1,
    \end{equation}
    and
    \begin{equation}
        \label{eq:bounds-proof-2}
        \lagp[1](x) = 1 - \frac{1}{\Delta^2} \sum_{j = 1}^{m - n - 1} \!\! x_j^2 - \frac{1}{\Delta} \sum_{j = m - n}^n \!\! x_j \le 1 -  \frac{1}{\Delta} \sum_{j = m - n}^n \!\! x_j \le 2n-m+2.
    \end{equation}
    Thus
    \begin{equation*}
        \Lambda_{\infty} = \max_{\norm{x}_{\infty} \le \Delta} @@ \abs[\big]{\lagp[1](x)} \le \max @ \set{n - 1, @@ 2n - m + 2},
    \end{equation*}
    which completes the proof of~\eqref{eq:bounds}.

    Finally, remark that the right-hand side in~\eqref{eq:bounds-proof-1} is attained by the vector whose entries are all~$\Delta$.
    Moreover, the right-hand side in~\eqref{eq:bounds-proof-2} is attained by the vector whose first~$m - n - 1$ entries are zero and whose remaining entries are all~$-\Delta$.
    Hence,~$\Lambda_{\infty} = \max @ \set{n - 1, @@ 2n - m + 2}$.
\end{proof}

\subsubsection{Some special cases}

Theorem~\ref{thm:bounds} provides the value of~$\Lambda_{\infty}$, which is either~$n - 1$ or~$n$ since~$m \ge n + 2$.
There are more cases where we can evaluate~$\Lambda_p$, as we will detail in the following.
First, when~$1 \le p \le 2$, $\Lambda_p$ actually equals the lower bound in~\eqref{eq:bounds}.

\begin{proposition}
    \label{prop:lambda-p-2}
    For any~$m \in \set{n + 2, @@ n + 3, @@ \dots, @@ 2n + 1}$ and any~$p \in[1, 2]$, we have
    \begin{equation}
        \label{eq:lambda-p-2}
        \Lambda_p = 1 + (2n + 1 - m)^{\frac{p - 1}{p}}.
    \end{equation}
\end{proposition}

\begin{proof}
    According to Theorem~\ref{thm:bounds}, it suffices to show that the right-hand side of~\eqref{eq:lambda-p-2} is an upper bound for~$\Lambda_p$.
    We will prove this using Lemma~\ref{lem:lambda-p}.

    Consider any~$x \in \mathcal{B}_p(\Delta)$.
    Note that~$\norm{x}_2 \le\norm{x}_p \le \Delta $ for~$p\in[1, 2]$.
    According to Lemma~\ref{lem:lagp} and the H{\"{o}}lder inequality, we have
    \begin{align*}
        -\lagp[1](x)    &= -1 + \frac{1}{\Delta^2} \sum_{j = 1}^{m - n - 1} \!\! x_j^2 + \frac{1}{\Delta} \sum_{j = m - n}^n \!\! x_j\\
                        &\le -1 + \frac{1}{\Delta^2}\norm{x}_2^2 + \frac{1}{\Delta} (2n+1-m)^{\frac{p-1}{p}} @ \norm{x}_p\\
                        &\le (2n+1-m)^{\frac{p-1}{p}}.
    \end{align*}
    Similarly,
    \begin{equation*}
        \lagp[1](x) \le 1 + \frac{1}{\Delta}\sum_{j = m - n}^n \abs{x_j} \le 1 + \frac{1}{\Delta} (2n + 1 - m)^{\frac{p-1}{p}} @ \norm{x}_p \le 1 + (2n + 1 - m)^{\frac{p - 1}{p}}.
    \end{equation*}
    Therefore, invoking Lemma~\ref{lem:lambda-p}, we have
    \begin{equation*}
        \Lambda_p = \max_{\norm{x}_p \le \Delta} @@ \abs[\big]{\lagp[1](x)} \le 1 + (2n + 1 - m)^{\frac{p - 1}{p}},
    \end{equation*}
    which concludes the proof.
\end{proof}

We can also evaluate~$\Lambda_p$ under the default and natural setting~$m = 2n + 1$ recommended by Powell.
To this end, recall the following elementary fact.

\begin{lemma}
    \label{lem:max-norm-pq}
    For any~$p \in [1, \infty]$ and $q \in [1, \infty]$, we have
    \begin{equation}
        \label{eq:max-norm-pq}
        \max_{\norm{x}_p \le 1} @@ \norm{x}_q = \max @ \set{1, @@ n^{\frac{1}{q} - \frac{1}{p}}}.
    \end{equation}
\end{lemma}

\begin{proof}
    If~$p$ or~$q$ is infinity, then~\eqref{eq:max-norm-pq} is straightforward to verify.
    So, we assume that both of them are finite.

    Consider the case where~$p \le q$.
    For~$x \in \mathcal{B}_p(1)$, we have~$\norm{x}_q^q \le \norm{x}_p^p \le 1$, and this bound is attained at the first coordinate vector~$e_1 \in \mathcal{B}_p(1)$, so that~\eqref{eq:max-norm-pq} holds in this case.

    We now suppose that~$p > q$.
    Let~$\ones \in \R^n$ denote the vector with all entries being one,~$r = p/q$, and~$s = r / (r - 1) = p / (p - q)$.
    For~$x \in \mathcal{B}_p(1)$, define~$z = (\abs{x_1}^q, \abs{x_2}^q, \dots, \abs{x_n}^q)$.
    According to the H{\"{o}}lder inequality, we have
    \begin{equation*}
        \norm{x}_q  = \big(\ones^{\T} z \big)^{\frac{1}{q}} \le \big(\norm{\ones}_s \norm{z}_r \big)
        ^{\frac{1}{q}} = n^{\frac{1}{sq}} \norm{x}_p \le n^{\frac{p - q}{qp}}.
    \end{equation*}
    Moreover, this bound is attained at~$x^\ast = n^{-\frac{1}{p}} \ones$, which proves~\eqref{eq:max-norm-pq} for~$p > q$.
\end{proof}

\begin{proposition}
    \label{prop:lambda-p-opt}
    For any~$p \in [1,\infty]$, if~$m = 2n + 1$, then
    \begin{equation*}
        \Lambda_p = \max @ \set[\big]{1, n^{\frac{p - 2}{p}} - 1}.
    \end{equation*}
\end{proposition}

\begin{proof}
    Since~$m = 2n + 1$, we have~$L_1(x) = 1 - \Delta^{-2}\norm{x}_2^2$ according to Lemma~\ref{lem:lagp}.
    Then, it is clear that
    \begin{equation}
        \label{eq:lambda-p-inf-p-1}
        \max_{\norm{x}_p \le \Delta} @ \lagp[1](x) = \max_{\norm{x}_p \le \Delta} @ \big( 1 - \Delta^{-2}{\norm{x}_2^2} \big) = 1.
    \end{equation}
    Moreover, according to Lemma~\ref{lem:max-norm-pq}, we have
    \begin{equation}
        \label{eq:lambda-p-inf-p-2}
        \max_{\norm{x}_p \le \Delta} @ -\lagp[1](x) = \max_{\norm{x}_p \le \Delta} @ \big( \Delta^{-2} \norm{x}_2^2 - 1 \big) = \max_{\norm{z}_p \le 1} @@ \norm{z}_2^2 - 1 = \max @ \set{0, @@ n^{\frac{p - 2}{p}} - 1}.
    \end{equation}
    The desired result is obtained by combining~\eqref{eq:lambda-p-inf-p-1} and~\eqref{eq:lambda-p-inf-p-2} with Lemma~\ref{lem:lambda-p}.
\end{proof}

Now we can show that~$\xpt[0]_{2n + 1}$ attains the minimal constant of well-poisedness in~$\mathcal{B}_p(\Delta)$ among all interpolation sets that contain the origin, provided that~$1 \le p \le 2$.
In this sense,~$\xpt[0]_{2n + 1}$ is an optimal interpolation set in~$\mathcal{B}_p(\Delta)$ for~$p \in [1, 2]$.
Indeed, the upper bound for~$p$ can be slightly larger than~$2$, as is detailed in Theorem~\ref{thm:optimset}.
Recall that $p = 2$ is of interest for most trust-region \gls{dfo} methods, including those by Powell.

\begin{theorem}
    \label{thm:optimset}
    Assume that~$m = 2n + 1$, and that either~$n \le 2$ or
    \begin{equation}
        \label{eq:pineq}
        1 \le p \le \frac{2 \log n}{\log (n/2)}.
    \end{equation}
    If an interpolation set containing~$0$ is~$\Lambda$-poised in~$\mathcal{B}_p(\Delta)$, then~$\Lambda \ge \Lambda_p$.
\end{theorem}

\begin{proof}
    Under the assumptions, we have~$\Lambda_p = 1$ according to Proposition~\ref{prop:lambda-p-opt}.
    If an interpolation set containing~$0$ is~$\Lambda$-poised in~$\mathcal{B}_p(\Delta)$, then we
    have~$\Lambda \ge 1$, because the Lagrange polynomial corresponding to~$0$ takes the value~$1$ at~$0$.
    Thus, the theorem holds.
\end{proof}
It is worth mentioning that Proposition~\ref{prop:lambda-p-2} also implies that~$\xpt[0]_{2n + 1}$ renders~$\Lambda_p = 1$ for~$p \in [1, 2]$.
Theorem~\ref{thm:optimset} shows that this equality holds for a larger range of~$p$, because the right-hand side in~\eqref{eq:pineq} is greater than~$2$ when~$n > 2$.

\subsubsection{Remarks and open questions}

Note that~$\Lambda_p$ can be regarded as a function of~$m$.
Theorem~\ref{thm:optimset} implies that~$m^\ast = 2n + 1$ minimizes this function if $p$ satisfies inequality~\eqref{eq:pineq}.
Thus, it is natural to ask whether~$2n+1$ minimizes~$\Lambda_p$ for any~$p \ge 1$, to which we do not have an answer yet.

The definition of~$\xpt[0]_m$ in~\eqref{eq:set-def} and~\eqref{eq:set-def-m} assumes that~$m \le 2n + 1$.
Even though larger values of~$m$ are not recommended in practice, Powell~\cite{Powell_2006} proposed an extension of~$\xpt[0]_m$ for~$m > 2n + 1$.
With such an extension, we can define~$\Lambda_p$ by~\eqref{eq:Lambdap} for any~$m \in \set{n + 2, @@ n +  3, @@ \dots, @@ (n + 1)(n + 2) / 2}$.
It is interesting to ask whether~$2n + 1$ still minimizes~$\Lambda_p$ for any~$p \ge 1$ after the extension.
We expect that the analysis will be more challenging than what we have done.
One of the challenges is that Lemma~\ref{lem:lambda-p} does not hold when~$m > 2n + 1$, and hence, the estimations of~$\Lambda_p$ will become more involved.

Another interesting open question is whether an interpolation set containing the origin can be~$\Lambda$-poised in~$\mathcal{B}_p(\Delta)$ with~$\Lambda < \Lambda_p$.
Theorem~\ref{thm:optimset} provides a negative answer when~$p$ satisfies~\eqref{eq:pineq}.
If the answer is negative for all~$p \ge 1$, then the optimality of~$2n + 1$ mentioned in the above two
paragraphs is true, and~$\xpt[0]_{2n + 1}$ is an optimal interpolation set in~$\mathcal{B}_p(\Delta)$ for all~$p \ge 1$.

\section{Conclusion}
\label{sec:conclusion}

We have analyzed the well-poisedness (in the minimum Frobenius norm sense) of an interpolation set that appears in many trust-region \gls{dfo} methods, particularly those by Powell~\cite{Powell_2006,Powell_2009}.
It is proved to be an optimal interpolation set under the default setting, because it minimizes the constant of well-poisedness in a ball centred at the starting point.
Our analysis justifies the natural configuration of this set from the viewpoint of interpolation theory.

The spirit of our analysis is similar to that of~\cite{Dodangeh_Vicente_Zhang_2016}, which proves that a widely used direction set is indeed optimal for directional direct search methods based on sufficient decrease.
While proving nothing surprising, this kind of investigation deepens our understanding of certain algorithmic strategies that we often employ but rarely ask why.

\paragraph*{Acknowledgement}
This paper corresponds to Section~2.5 of the PhD thesis of Tom M. Ragonneau~\cite{Ragonneau_2022}, co-supervised by Zaikun Zhang and Professor Xiaojun Chen from The Hong Kong Polytechnic University.
Both authors are very grateful to Professor Chen for her support, encouragement, and guidance during the thesis.
Zaikun Zhang would like to thank the late Professor Oleg Burdakov for his friendship.

\paragraph*{Disclosure statement}
The authors report that there are no competing interests to declare.

\paragraph*{Funding}
This work was funded by the University Grants Committee of Hong Kong under projects PF18-24698 (Hong Kong PhD Fellowship Scheme), PolyU 253012/17P, PolyU 153054/20P, PolyU 153066/21P, and PolyU 153086/23P.
It was also supported by The Hong Kong Polytechnic University under projects P0009767, P0038928, P0045598, and the CAS-Croucher Funding Scheme for \textquote[][.]{CAS AMSS-PolyU Joint Laboratory of Applied Mathematics: Nonlinear Optimization Theory, Algorithms and Applications}

\bibliography{optimset}

\end{document}